\mag\magstep1

\magnification1200
\input amssym.def 
\input amssym.tex 
\def\SetAuthorHead#1{}
\def\SetTitleHead#1{}
\def\NoindentAfter{\everypar={\setbox0=\lastbox\everypar={}}}
\def\H#1\par#2\par{{\baselineskip=15pt\parindent=0pt\parskip=0pt
 \leftskip= 0pt plus.2\hsize\rightskip=0pt plus.2\hsize
 \bf#1\unskip\break\vskip 4pt\rm#2\unskip\break\hrule
 \vskip40pt plus4pt minus4pt}\NoindentAfter}
\def\HH#1\par{{\bigbreak\noindent\bf#1\medbreak}\NoindentAfter}
\def\HHH#1\par{{\bigbreak\noindent\bf#1\unskip.\kern.4em}}
\def\th#1\par{\medbreak\noindent{\bf#1\unskip.\kern.4em}\it}
\def\endth{\medbreak\rm}
\def\pf#1\par{\medbreak\noindent{\it#1\unskip.\kern.4em}}
\def\df#1\par{\medbreak\noindent{\it#1\unskip.\kern.4em}}

\let\Roster\bgroup\let\endRoster\egroup
\def\\{}\def\text#1{\hbox{\rm #1}}

\def\MaxReferenceTag#1{}
\def\qedbox{\vrule width2mm height2mm\hglue1mm\relax}
\def\qed{\ifmmode\qedbox\else\hglue5mm\unskip\hfill\qedbox\medbreak\fi\rm}

\def\cite#1{{\bf[#1]}}
\def\Em#1{{\it #1\/}}
\def\Bib#1\par{\bigbreak\bgroup\centerline{#1}\medbreak\parindent30pt
 \parskip2pt\frenchspacing\par}
\def\endBib{\par\egroup}
\newdimen\Overhang
\def\rf#1{\par\noindent\hangafter1\hangindent=\parindent
     \setbox0=\hbox{[#1]}\Overhang\wd0\advance\Overhang.4em\relax
     \ifdim\Overhang>\hangindent\else\Overhang\hangindent\fi
     \hbox to \Overhang{\box0\hss}\ignorespaces}

\def\Coordinates{\bigbreak\bgroup\parindent=0pt\obeylines}
\def\endCoordinates{\egroup}

\let\eval\Overline
\SetTitleHead{}
\SetAuthorHead{}

\def\Title {A Note on Context[C Sensitive Languages\\ and Word Problems}

\def\Author{Michael Shapiro\footnote{*}{I wish to thank the NSF for support.}}

\H \Title

\Author

\HH

In \cite{AS}, Anisimov and Seifert show that a group has a regular
word problem if and only if it is finite.  Muller and Schupp \cite{MS}
(together with Dunwoody's accessibility result \cite{D}) show that a
group has context free word problem if and only if it is virtually
free.  In this note, we exhibit a class of groups where the word
problem is as close as possible to being a context sensitive language.
This class includes the automatic groups of \cite{ECHLPT} and is
closed under passing to finitely generated subgroups.  Consequently,
it is quite large.  For example, it contains all finitely generated
subgroups of the $n$-fold product of free groups, $F_2 \times
\dots\times F_2$.  For $n=2$, these include groups which are not
finitely presented, and for $n>2$, these include groups which are
$FP_n$ but not $FP_{n+1}$.

Let us make clear what we mean by saying that the word problem is as
close as possible to being a context sensitive language.  Recall that
a context sensitive language cannot contain the empty word $e$.  Since
the empty word is always an element of the word problem, strictly
speaking, the word problem can never be a context sensitive language.
So we will abuse terminology and say that the word problem is context
sensitive if, after deleting the empty word, it is context sensitive.
We feel that this is not a grievous abuse: in any practical situation
where one is trying to either decide or enumerate the word problem,
the empty word is the least of one's problems!

There are two ingredients to our Theorem.  One is the notion of a
\Em{asynchronous combing} of a group (see below). The other is the
following characterization of context sensitive languages.  Given a
language $L$, $L-\{e\}$ is context sensitive if and only if $L$ is the
language of a nondeterministic linear bounded Turing machine\footnote
{\dag}{Indeed, as Neumann has pointed out to me, one might wish to
cure this mismatch by the following method.  We could change the
definition of a context sensitive language to allow a finite number of
start words rather than a single start symbol.  If none of these is
$e$, each could be produced from a single start symbol by a single
production rule.  Thus the result of this change of definition would
be to include $L$ if and only if $L-\{e\}$ is context sensitive under
the old definition.}.  (See, for example, \cite{HU}.) Thus, to see
that the word problem is context sensitive, we must exhibit an
algorithm which, given $w$ with $\eval w=1$, verifies membership in
the word problem using an amount of space which is linear in the
length of $w$.  In fact, the process we will describe gives a
deterministic linear bounded automaton in the case where the combing
language is also the language of a deterministic linear bounded
automaton.  In this case, our algorithm acts to decide the word
problem rather than merely verify membership in the word problem.

We start by fixing our terminology.  Given a group $G$ and a finite
monoid generating set $A$, we take $A^*$ to be the free monoid on $A$.
For each $w \in A^*$ we denote the length of $w$ by $\ell(w)$.  The
\Em{empty word} is the unique word of length $0$ and we denote it by
$e$. We map $A^*$ to $G$ by the monoid homomorphism which takes each
letter of $A$ to its value in $G$.  We denote this map by $w \mapsto
\eval w$.  We call $\{ w \in A^* \mid \eval w = 1 \}$ the \Em {word
problem}.  We will assume that $A$ is supplied with an involution
denoted by $a \mapsto a^{-1}$ and that $\eval{a^{-1}}=(\eval a)^{-1}$
for all $a \in A$.  This allows one to build the \Em{Cayley graph}
$\Gamma$ of $G$ with respect to $A$.  This is the labelled directed
graph whose vertices are the elements of $G$ and whose edges are
$\{(g,a,g') \mid g,g' \in G, a \in A, g'=g \eval a\}$.  Each edge
$(g,a,g')$ is labelled by $a$.  Elements of $A^*$ are called words,
and each word now labels a unique edge path of $\Gamma$ based at $1
\in \Gamma$.  Declaring each edge isometric to the unit interval
induces the \Em{word metric} $d(\cdot,\cdot)$ on $G$ and a \Em{length
function} $\ell(g) = d(1,g)$.  We call a subset of $A^*$ a
\Em{language}.  We call a language $L$ a \Em{normal form} if $\eval L
= G$. (We do not demand that this is a bijection.)  We call a normal
from $L$ an \Em{asynchronous combing} if there is a constant $K$ so
that for any $w,w' \in L$ with $d(\eval w, \eval{w'}) \le 1$, we can
find monotone reparameterizations of $[0,\infty)$ $t \mapsto t'$ and
$t \mapsto t''$ so that for all $t$, $d(w(t'),(w'')) \le K$.  We say
that  $D$ is a \Em{departure function} for $L$ if for any $w=xyz \in
L$, $\ell(\eval y) \ge n$ whenever $\ell(y) \ge D(n)$.  We say a
language $L$ is \Em{short} if there are $\lambda$ and $\epsilon$ so
that if $w \in L$ then $\ell(w) \le \lambda \ell(\eval w) + \epsilon$.
In addition, we say that $L$ consists of
$(\lambda,\epsilon)$\Em{-quasigeodesics} if for any $w=xyz \in L$,
$\ell(y) \le \lambda \ell(\eval y) + \epsilon$.

\th Theorem

Suppose that $H$ is a finitely generated subgroup of $G$ and suppose
that $G$ possesses a short asynchronous  context sensitive combing
with a departure function.  Then $H$ has a context sensitive word
problem. \endth

\th Corollary

A finitely generated subgroup of an automatic group has context
sensitive word problem. \endth

Indeed, we may replace ``automatic'' by the less popular but equally
serviceable class ``quasigeodesic asynchronously automatic'' \cite{N}.

\pf Proof

As we remarked above, we need to give a linear bounded algorithm for
verifying membership in the word problem.  Thus, it suffices to see
that $G$ has context sensitive word problem.  For suppose that $H$ is
generated by $B=\{h_1, \ldots , h_k\}$.  We choose $w_1, \ldots , w_k
\in A^*$ so that $\eval {w_i}=h_i$ for $i-1, \ldots , k$.  Then, given
$w'=h_{i_1} \ldots h_{i_n} \in B^*$, we replace this by
$w=w_{i_1}\ldots w_{i_n} \in A^*$.  This has increased length by at
most a factor of $\max\{\ell(w_i)\}$.  We now appeal our linear
bounded algorithm to determine if $\eval w =1$, and this will be
linearly bounded in the length of our original word $w'$. 

Let $L$ be our combing, and  suppose that we are given $w=a_1 \ldots
a_n \in A^*$.  Suppose that for $i=1, \ldots , n$, $u_i \in L$ and
$\eval {u_i}=w(i)$.  Then for each $i$, $\ell(u_i) \le \lambda n +
\epsilon$, where $\lambda$ and $\epsilon$ are the constants which
assure us that $L$ is short.  Further, $\eval w = 1$ if and only if
$\eval {u_n} =1$ and this happens if and only if $u_n$ is one of
finitely many words (all of length at most $\epsilon$.)  So once we
have found $u_n$, it is easy to determine if $\eval w=1$.  Thus it
suffices to see that we can find $u_n$ in a linearly bounded manner,
and to do this, it will suffice to show that we can find each
$u_{i+1}$ from $u_i$ and $a_{i+1}$ in a manner which is linearly
bounded in terms of $\ell(u_i)$, since this latter is itself linearly
bounded in terms of $\ell(w)$.

To do this, we start enumerating the words of $A^*$, say in short-lex
order, and test each one to see if it is an element of $L$.  Since $L$
is context sensitive, we can do this in a linearly bounded fashion.
When we find a word $u$ in $L$, we must check to see if it can be
taken as $u_{i+1}$.  That is, we must check whether $\eval u = w(i+1)
= \eval {u_i a_{i+1}}$.

Now $L$ is an asynchronous combing with a departure function.  Thus
for each $a \in A$, one can build an asynchronous two tape finite
state automaton which determines when given $u, u' \in L$ whether or
not $\eval u = \eval{u'a}$.  (For details see \cite{ECHLPT} or
\cite{BGSS}.)  Thus the decision as to whether or not to take $u$ as
$u_{i+1}$ can be made using an  amount of memory which is bounded by
a global constant.

If we have not found $u_{i+1}$, we go on to the next element of $A^*$
and discard $u$.  Eventually, we find $u_{i+1}$, and we need never
check any word of length longer than $\lambda (\ell(u_i)+1)+\epsilon$.
Since no $u_i$ has length longer than $\lambda n +\epsilon$ we shall
eventually find $u_n$ in a linearly bounded way.  \qed

The Corollary follows by noting that an automatic group has an
automatic structure with uniqueness.   This will consist of
quasigeodesics and hence is short and has a departure function.  We
can assume it does not contain $e$.  This, together with the fact
that it is regular, ensures that it is context sensitive.

\medskip
In the case where $L$ is an automatic structure with uniqueness,
\cite{ECHLPT} show that $u_{i+1}$ can be found from $u_i$ by a process
whose time is linearly bounded in $\ell(u_i)$.  This gives a method
for solving the word problem in quadratic time.  In a similar vein, if
$G$ is a direct product of word hyperbolic groups, the word problem
for $H$ can be solved in linear time using pushdown automata \cite{S}.

\Bib{References}
\MaxReferenceTag{ECHLPT}

\rf{AS} A.V.~Anisimov and F.D.~Seifert, Zur algebraischen
charateristik der durch kontextfreie Sprachen definierten Gruppen,
Elektron.\ Informationsverarb.\ Kybernet. 11 (1975) 695---702.

\rf{BGSS} G.~Baumslag, S.M.~Gersten, M.~Shapiro and H.~Short,
Automatic groups and amalgams, Journal of Pure and Applied Algebra 76
(1991), 229---316.

\rf{D} W.~Dunwoody, The Accessibility of finitely presented groups,
Inventiones Mathematica, (1985), 449---457.

\rf{ECHLPT} D.B.A.~Epstein, J.W,~Cannon, D.F.~Holt, S.V.F.~Levi,
M.S.~Paterson and W.P. Thurston, ``Word Processing in Groups,'' Jones
and Bartlett Publishers, Boston, 1992.

\rf{HU} J.E.~Hopcroft and J.D.~Ullman, ``Introduction to Automata
Theory, Languages, and Computation," Addison Wesley, Reading, 1979.

\rf{MS} D.E.~Muller and P.E.~Schupp, Groups, the theory of ends, and
context-free languages, Journal of Computer and System Sciences, 26
(1983), 295---310.

\rf{N} W.D.~Neumann, Asynchronous combings of groups, International
Journal of Algebra and Computation, 2 (1992), 179---185.

\rf{S} H.~Short {\it et al.},  Notes on word hyperbolic groups, in
Group Theory from a Geometric Viewpoint, E.~Ghys, A.~Haefliger,
A.~Verjovsky eds., World Scientific, 1991.

\bye